\documentclass[12pt]{article}

\usepackage[english]{babel}
\usepackage[utf8]{inputenc}

\usepackage[intlimits]{amsmath} 
\usepackage{amsthm,mathrsfs}    
\usepackage[comma,authoryear]{natbib}
\usepackage{url,color}
\numberwithin{equation}{section}
\usepackage{amssymb,amsmath}
\usepackage{subfigure}
\usepackage{bm}
\usepackage{booktabs} 
\usepackage{graphicx,graphics,epsfig}
\usepackage{float}

\theoremstyle{plain}
\newtheorem{thm}{Theorem}[section]

\newcommand{\bthm}{\begin{thm}}
\newcommand{\ethm}{\end{thm}}

\newcommand{\bpf}{\begin{proof}}
\newcommand{\epf}{\end{proof}}

\theoremstyle{definition}

\newcommand{\GF}{GraField}
\newcommand{\ER}{Erd{\H{o}}s-R{\'e}nyi\,\,}

\usepackage[a4paper,centering,text={5.95in,9.25in}]{geometry}

\usepackage[compact]{titlesec}
\newcommand{\bib}{\bibliography{ref-bib}\bibliographystyle{Chicago}}
\usepackage{deep-macro}

\begin{document}
\begin{center}
{\Large{\bf Strength of Connections in a Random Graph:\\\vskip.25em Definition, Characterization, and Estimation}} \\[.2in] 
Subhadeep Mukhopadhyay\\
Department of Statistics, Temple University \\ Philadelphia, Pennsylvania, 19122, U.S.A.\\
\texttt{deep@temple.edu}\\
\vskip3em
{\bf ABSTRACT}\\
\end{center}
\vspace{-1em}
How can the `\emph{affinity}' or `\emph{strength}' of ties of a random graph be characterized and compactly represented? How can concepts like \emph{Fourier and inverse-Fourier} like transform be developed for graph data? To do so, we introduce a new graph-theoretic function called `Graph Correlation Density Field' (or in short \texttt{GraField}), which differs from the traditional edge probability density-based approaches, to completely characterize tie-strength between graph nodes. Our approach further allows frequency domain analysis, applicable for both directed and undirected random graphs.

\vspace*{.08in}
\vspace*{.08in}


\section{Introduction}
\cite{lovasz2006} introduced the concept of Graphon, as a functional modeling tool to compactly describe the edge probabilities of a random graph. Characterization and estimation based on graphon is now a fast-growing branch of (nonparametric) statistics and graph theory; see for instance \cite{airoldi2013}, \cite{olhede2013}, \cite{chatterjee2014}, \cite{frieze1999}, and \cite{lovasz2012}.

However, the fundamental quantity of interest in graphs or networks is \emph{not} the `edge probabilities,' rather it is the connections (or interaction) between pairs of vertices. In the literature of random graph modeling there \emph{has not been any attempt} to define such `strength' of connections (or interaction) between pairs of vertices. In this paper, our aim is to construct a function that characterizes this \emph{connectedness or tie-strength} (instead of of representing the probability of a tie or edge between any two nodes) in a compact way and to study its fundamental role in random graph modeling. The central goal of this paper is to answer:
\vskip.15em
\emph{How can such implicit affinity function be constructed and estimated from explicitly observed links?}
\vskip.12em
We introduce a fundamental object, called `\underline{Gra}ph Correlation Density \underline{Field}' (\texttt{\GF}), as a tool to provide a compact functional description of graph affinity (how related any two vertices are), formulated precisely in Section 2.3 and 2.4. Orthogonal decomposition of this function is also derived to efficiently represent and nonparametrically estimate the underlying correlation field. We also address the problem of graph \emph{compression} and \emph{sparse} representation learning to model the interactions between large numbers of nodes - a major research challenge in modern graph theory.

Furthermore, we show that how the concept of \texttt{\GF} provides a systematic and comprehensive framework for frequency-domain analysis of random graphs, which is \emph{analogous} to its time series counterparts. In particular, we are interested in the following question, which plays a fundamental role in graph signal processing:
\vskip.15em
\emph{How can concepts like of Fourier and inverse-Fourier like transform be developed for graph data?}
\vskip.15em

Our nonparametric approximation scheme of \texttt{\GF} hinges upon a specially designed orthogonal transform basis. Section 2.1 introduces the theory and construction of this new kind of discrete (piecewise constant) orthonormal basis set called LP Graph basis. We define the concept of LP graph transforms (analogue of Fourier transform in time series) in Section 2.2 to allow a frequency interpretation of graph data parallel to time series analysis that captures graph connectivity patterns and plays a vital role in compression. We also show that the LP graph transform coefficients can be interpreted as discrete orthogonal transforms of the fundamental Graph Correlation Density function.

We apply this new graph correlation field theory to solve two graph modeling problems in Section 3. The first step of modeling is to identify the null model. In time series analysis, the `flat' shape of power spectral density and the plot auto-correlations provide an elegant diagnostic tool for white noise process (which is the null model in the context of time series). Section 3.1 addresses the following question:
\vskip.15em
\emph{Can we develop an analogous portmanteau graphical diagnostic test (similar in spirit to time series) for null random graph models?}
\vskip.15em
We discuss two graphical diagnostic tests that resemble the sample auto-correlation function and spectral density-based approach for white noise testing in times series context. We show that for a null graph model, the shape of `Graph Correlation Density' is `flat' and the sample LP-graph transforms behave asymptotically normal distribution with mean $0$ and variance $1/N$, remarkably similar to its time series counterparts  \citep{brockwell2009}.

Section 3.2 explores the connection between two fundamental functions of modern graph theory: Graph Correlation Density and the traditional edge probability-based concept, Graphon \citep{lovasz2012}. As a practical utility of this connection we answer the following question:
\vskip.15em
\emph{Can we develop a sparse estimation strategy that will provide `smooth' nonparametric estimate of Graphon?}
\vskip.15em
Our research advances the current state-of-the-art nonparametric graphon modeling (that yields histogram-based approximation) by offering a low-rank efficient representation-based smoothing algorithm. We end with final conclusions and future research directions in Section 4.
\section{Theory and Estimation}
\subsection{Discrete Basis Construction}
\begin{figure}[thb]
 \centering
 \vspace{-.2em}
 \includegraphics[height=\textheight,width=.49\textwidth,keepaspectratio,trim=1cm 1cm .4cm 1cm]{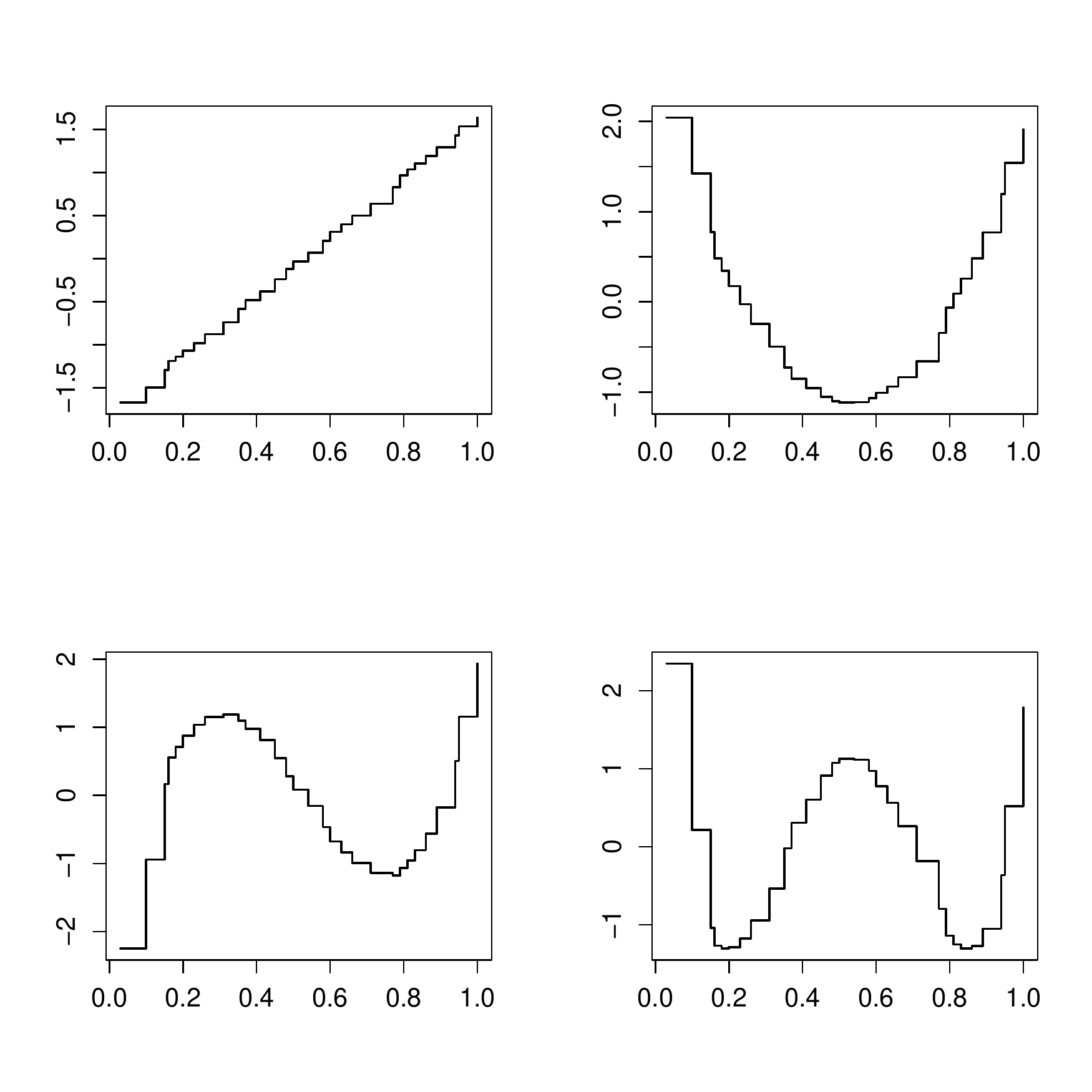}
 \includegraphics[height=\textheight,width=.49\textwidth,keepaspectratio,trim=.4cm 1cm 1cm 1.2cm]{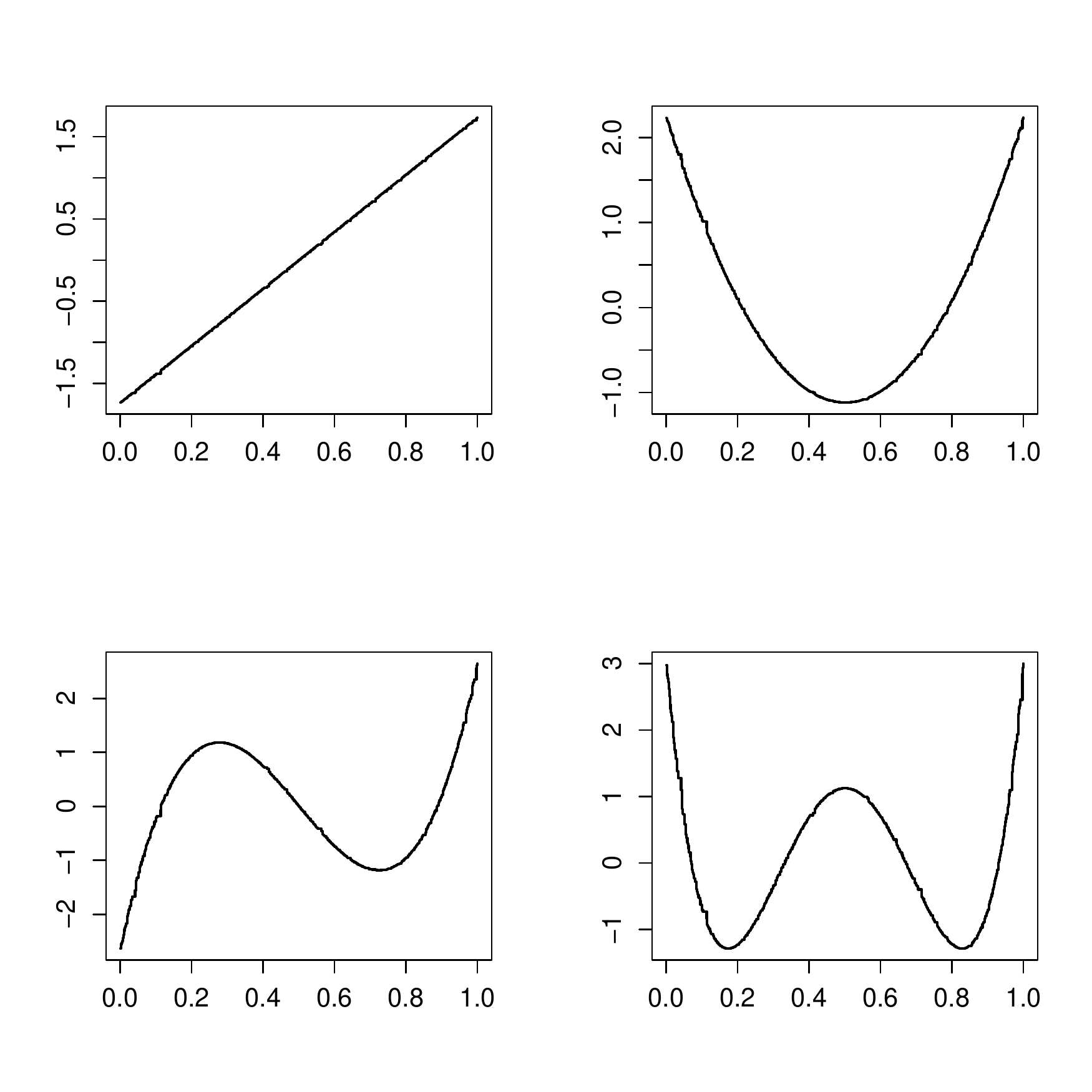}
 \caption{(a) Left: Piecewise constant shapes of first four LP orthonormal basis for the bipartite graph with $n=30$ (divided equally into two disjoint subsets) and $p=.25$; (b) Right: first four basis functions for political blog data that resemble orthonormal shifted Legendre polynomials on $(0,1)$.}
  \label{fig:score}
\end{figure}
The traditional approach constructs nonparametric orthonormal graph basis by performing singular-value decomposition (SVD) of Laplacian or Adjacency matrix (both matrices are of size $n \times n$, where $n$ is the graph size). The question that naturally arises is as follows:
\vskip.15em
\emph{Can we bypass this computationally expensive SVD of a (possibly) large matrix to find the graph Fourier-like discrete orthonormal bases}?
\vskip.15em
In what follows, we propose an algorithm that finds nonlinear adaptive discrete graph basis avoiding SVD altogether - computationally fast. Our basis is nonparametric (\emph{not fixed} like discrete Fourier basis, wavelets, or curvelets, etc.) in the sense that its shape depends on (or adapts to) the given graph structure. Furthermore, these custom design basis functions will act as an expansion basis for our functional graph object $\cd(u,v;\cG)$ `graph correlation density'.

Here, we provide the recipe for how to construct discrete piecewise constant orthonormal basis functions, specially designed for the given random graph. Let us start with the construction of orthonormal polynomials associated with an arbitrary discrete random variable $X$. For a given discrete random variable $X$, we construct the basis functions by Gram Schmidt orthonormalization of powers of $T_1(x;X)$, given by
\beq T_1(x;X) \,=\,\dfrac{\sqrt{12} [\Fm(x;X) - .5]}{\sqrt{1-\sum_x p^3(x;X)}}, \eeq
where
\beq
\Fm(x;X)=F(x;X)-.5p(x;X),\, p(x;X)=\Pr[X=x], \,F(x;X)=\Pr[X \leq x].
\eeq
By construction, these basis functions satisfy the following orthonormality conditions:
\beq \label{eq:ortho}
\sum_i T_j(x_i;X) p(x_i;X)=0, \sum_i T^2_j(x_i;X) p(x_i;X)=1, \sum_i T_j(x_i;X)  T_k(x_i;X) p(x_i;X) = 0.
\eeq
We define LP orthonormal basis in the \emph{unit interval} $[0,1]$ (rescaling operation) by evaluating $T_j$ at $Q(u;X)$, quantile function of the random variable $X$
\beq S_j(u;X)\,=\,T_j(Q(u;X);X), ~0<u<1. \eeq

For a given directed (and possibly  weighted) random graph of size $n$, we construct two sets of nonlinear adaptive discrete orthogonal graph bases $\{S_1(u;X),S_2(u;X), \ldots S_{n-1}(u;X)\}$ and $\{ S_1(v;Y),S_2(v;Y), \ldots S_{n-1}(v;Y) \}$ based on in-degree and out-degree sequences
\beq
p(x;X)=\sum_{y=1}^n A(x,y)/N , p(y;Y)=\sum_{x=1}^n A(x,y)/N,  p(x,y;\cG)=A(x,y)/N, N=\sum_{x,y=1}^n A(x,y),
\eeq
where $A$ is the  (weighted) adjacency matrix of the graph $\cG$. The Eq \eqref{eq:ortho} guarantees that both basis sets are automatically \emph{degree-weighted}. Note that for undirected graphs we only have one set of basis due to the symmetry.

Figure \ref{fig:score}(a) shows the `piecewise constant' (nonlinear) shapes of the first four orthonormal basis functions of the bipartite graph with $n=30$ (divided equally into two disjoint subsets) and connection probability $p=.25$. However, the LP graph bases exhibit a `smooth continuous shape' for the  political blog data \citep{adamic2005}, as shown in Figure \ref{fig:score}(b). This is a directed graph with $n=1490$ and the total of directed edges (incoming and outgoing links) $19,090$. This data was collected immediately after the $2004$ presidential election to understand the connection patterns (as indicated by hyperlinks) between two political web-blogs (liberal and conservative).

One unique aspect of our theory is that it enables us to pass from the `piecewise-constant discrete' orthonormal basis of $\cR^n$ (for small graphs) to an `approximately smooth' continuous orthonormal basis of $L^2(0,1)$ (for large graphs). In the \emph{continuum limit} (as dimensional of the graph $n \rightarrow \infty$), our bases exhibit a `\emph{universality}' property in the sense that it approaches to shifted orthonormal \underline{L}egendre \underline{P}olynomials on $(0,1)$. To emphasize this universal limiting shape, we call it `LP' graph basis functions.
\subsection{Orthogonal Graph LP Transform}
We introduce a Fourier-like transform, called graph LP transform (valid for directed/undirected graphs) that provides a novel scheme for capturing the structural information of a graph using few coefficients, thus achieving the goal of graph compression.

Compute LP Graph transforms by taking cross-covariance (which can be interpreted as correlation due to \eqref{eq:ortho}, thus we can alternatively call it LP graph correlations) between the LP-basis functions $T_j(X;\cG)$ and $T_k(Y;\cG)$ for $j=\{1,2,\ldots,n-1\}$ and $k=\{1,2,\ldots,n-1\}$
\beq
\LP[j,k;\cG]~=~ \Ex[T_j(X;\cG)T_k(Y;\cG)]~=~\mathop{\sum\sum}_{x,y\in \{1,2,\ldots,n\}} p(x,y;\cG) T_j(x;\cG) T_k(y;\cG)~\equiv~\langle T_j, T_k \rangle_{L^2(p)}.
\eeq
We go from the original graph domain to the LP frequency domain since important structural information of a graph may be compressed into fewer discrete LP transform coefficients. Strictly speaking, the topology for structured graphs tends to be concentrated in a few low-frequency (low values of $j$ and $k$) LP graph transforms - leading to efficient representation. Each LP graph correlations $\LP[j,k;\cG]$ captures different connection patterns of the graph. As we will see in the next section, these LP transform coefficients can be interpreted as discrete orthogonal transforms of the graph correlation density kernel.

\subsection{Graph Correlation Density Field}
\begin{figure}[thb]
 \centering
 \vspace{-.2em}
 \includegraphics[height=\textheight,width=.49\textwidth,keepaspectratio,trim=1.8cm 1cm 1cm 1.8cm]{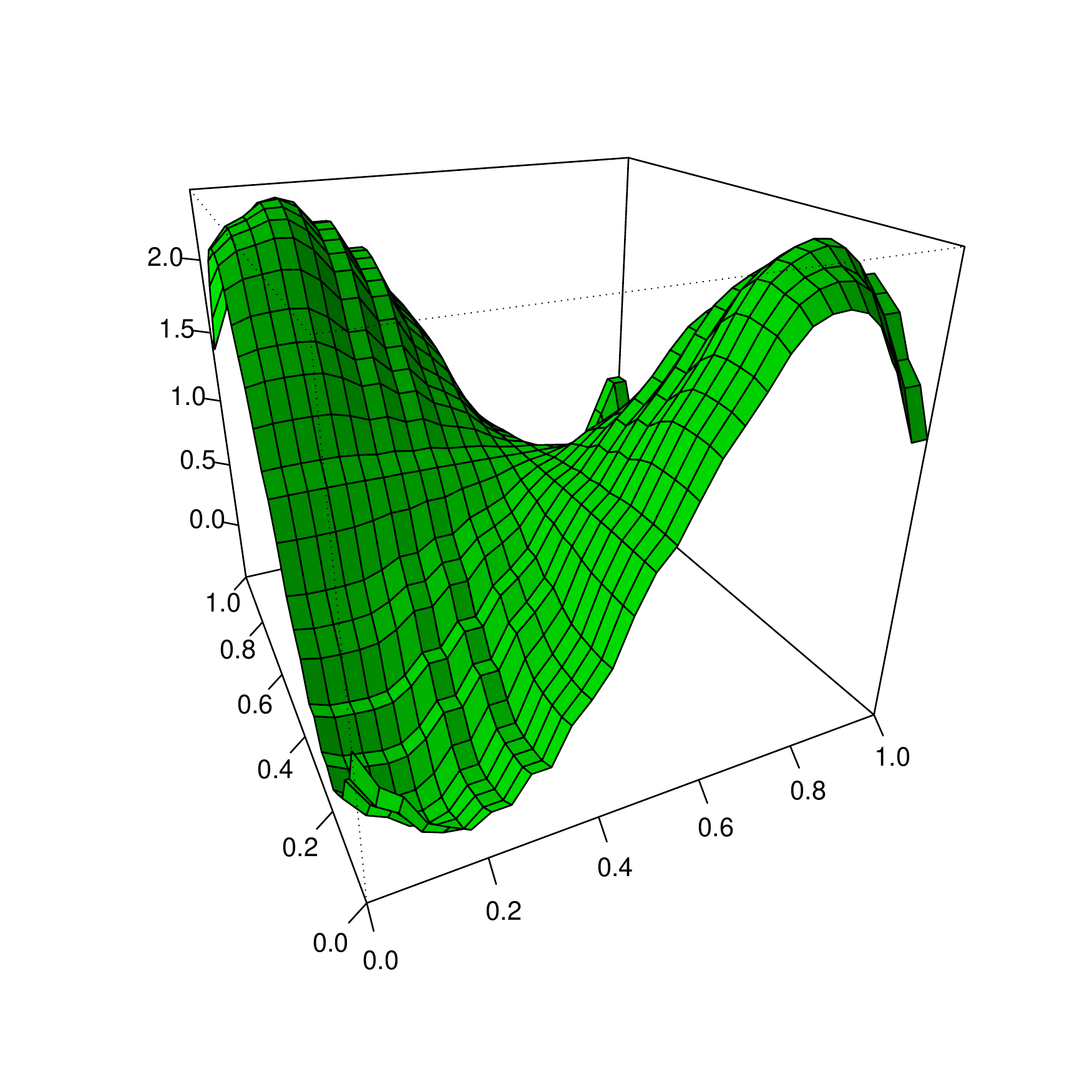}
 \includegraphics[height=\textheight,width=.49\textwidth,keepaspectratio,trim=1.5cm 1cm 1.5cm 1.5cm]{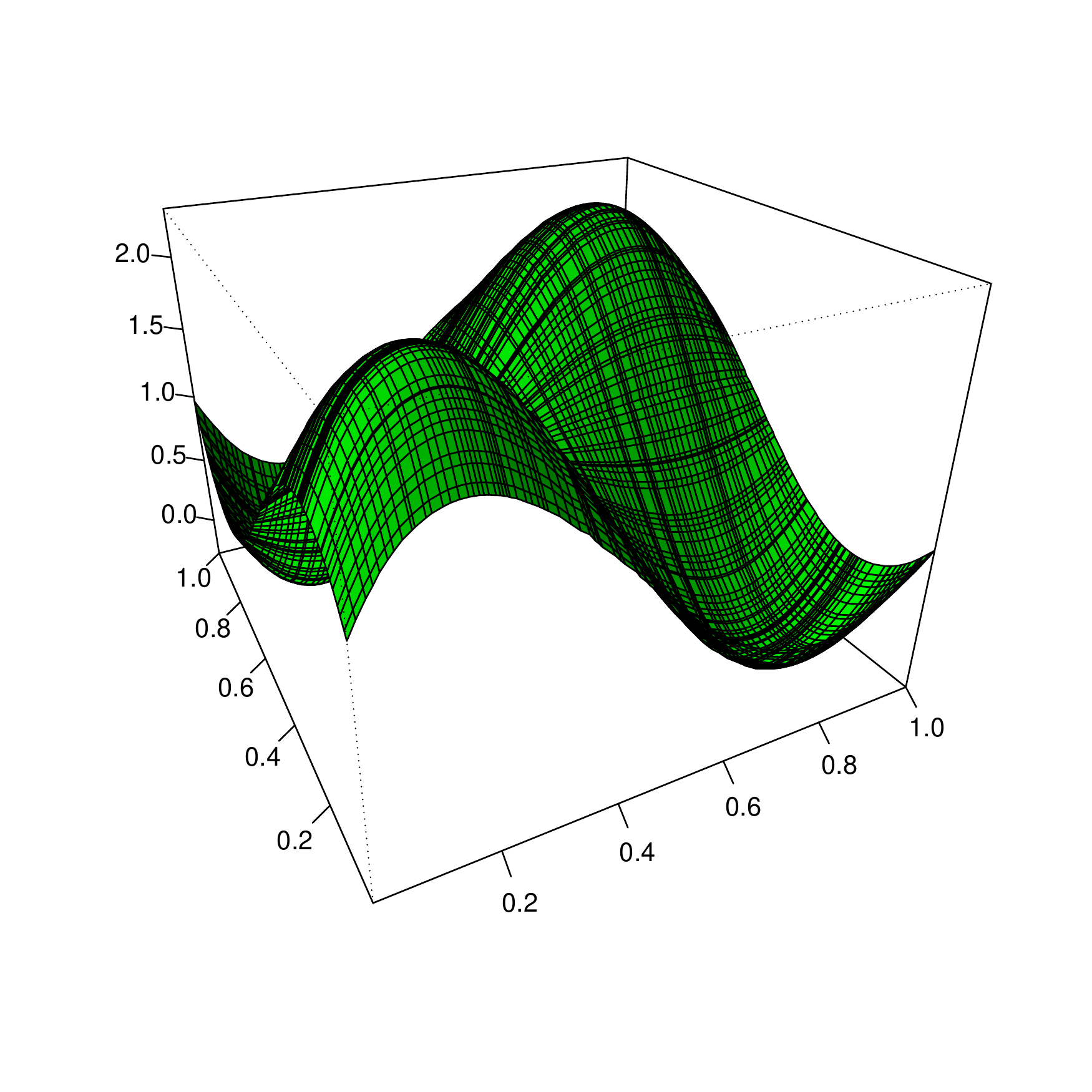}
 \caption{(color online) Nonparametric smooth graph correlation density estimate for (a) Left: bipartite graph  with $n=30$ and $p=.25$, and (b) Right:  directed political blog graph data with $n=1490$. Note the transition from bivariate step-function (small graph of size $n=30$) to smooth continuous shape (for large graph with $n=1490$ nodes) of the estimated \GF \,\,$\widehat \cd(u,v;\cG)$.}
  \label{fig:cop1}
\end{figure}
To completely characterize and compactly represent the `affinity' or `strength' of ties (or interaction) between every pair of vertices in the graph, we introduce a new functional statistical tool, called `Graph Correlation Density Field' (\texttt{\GF}), denoted by $\cd(u,v;\cG)$. We also discuss estimation of this (implicit) tie-strength function to capture the underlying interactions from the observed adjacency matrix $A$.
\begin{thm} \label{thm:dual} The two equivalent dual representations of orthogonal LP graph transform is given by
\beq
\LP[j,k;\cG]\,~~=\,~~\big\langle T_j, T_k \big \rangle_{L^2(p)}\,~~=\,~~\big\langle S_j, \int_0^1\cd\,S_k \big\rangle_{L^2[0,1]}~~~~(j,k=1,\ldots,n),
\eeq
where the square integrable function $\cd:[0,1]^2 \rightarrow \cR_+ \cup \{0\}$ is given by
\beq \label{eq:gcddef}
\cd(u,v;\cG)\,=\,\dfrac{p\big( Q(u;X), Q(v;Y);\cG \big)}{p\big( Q(u;X) \big) p\big(  Q(v;Y) \big)},
\eeq
where $u=F(x;X),v=F(y;Y)$ for $x,y \in \{1,2,\ldots,n\}$.
\end{thm}

The following fundamental representation theorem can be interpreted as \emph{inverse} LP \emph{transform} (in light of Theorem \ref{thm:dual}), which converts the set of LP-coefficients $\LP[j,k;\cG]$ into a density function over unit square.

\begin{thm} \label{thm:gcd} The function $\cd$ defined in Eq. \eqref{eq:gcddef} for a directed graph with size $n$, is a proper density function that admits the following orthogonal decomposition in terms of piecewise-constant product LP-basis over unit square for $u=F(x;X),v=F(y;Y)$ and $x,y \in  \{1,2,\ldots,n\}$.
\beq \label{eq:GCOR}
\cd(u,v;\cG)-1~=~\sum_{j=1}^{n-1}\sum_{k=1}^{n-1}   \LP[j,k;\cG]\, S_j(u;X)S_k(v;Y),
\eeq

\end{thm}

This can be considered as a frequency (LP) domain functional representation of graph topology, which computes the underlying (unobserved) `\emph{affinity} to connect' (not the un-normalized `raw' probability of an edge) between graph nodes. As the total area under this non-negative correlation function is one, we call it `density function'. LP space description of the graph structure \eqref{eq:GCOR} allows a Fourier-like interpretation. This also shows that the LP-coefficients can be used to reconstruct graph topology, captured by $\cd(u,v)$, which can be represented as the weighted sum of custom-designed product-LP graph basis functions.  Perform `smooth' estimation by selecting the `large' LP-components $\LP[j,k;,\cG]$ in a data-driven way (discussed in the next section).

We quantify the graph structure by measuring the departure from uniformity
\beq \label{eq:lpinfor}
\LPINFOR(\cG)\,=\,\sum_{j=1}^{n-1} \sum_{k=1}^{n-1} \Big|\LP[j,k;\cG]\,\Big|^2\,=\,\int_{[0,1]^2} \cd^2 - 1
\eeq
The equality follows directly from Parseval’s theorem.

Fig \ref{fig:cop1} shows the nonparametrically estimated graph correlation density $\widehat{\cd}(u,v;\cG)$ for both the bipartite graph ($n=30$) and political blog ($n=1490$) examples. The modes represent the higher connectedness. For the bipartite graph, it displays an underlying strong connection between the two disjoint sets and almost no interaction among the nodes from the same group. On the other, hand the web-blog data (directed network) clearly indicates the effect of political polarization. The different sharpness of the two modes indicate the contrasting connection strength within each of the political groups: liberal and conservative. Section 3.2 will further exemplify this point. In addition, note that for small graphs the corresponding shape of graph correlation density is bivariate non-negative (but not necessarily symmetric, as in the case of directed graphs) step-functions while the continuum view (for large graphs as $n \rightarrow \infty$) of $\cd(u,v;\cG)$  is a smooth continuous function over unit interval, hence we can also call it a ``correlation field''. Our theory provides a smooth macroscopic description of the microscopic interactions between the nodes (via adjacency matrix) in a compact way.
\begin{figure}[thb]
 \centering
 \includegraphics[height=.3\textheight,width=\textwidth,keepaspectratio,trim=1cm 1cm 1cm 1cm]{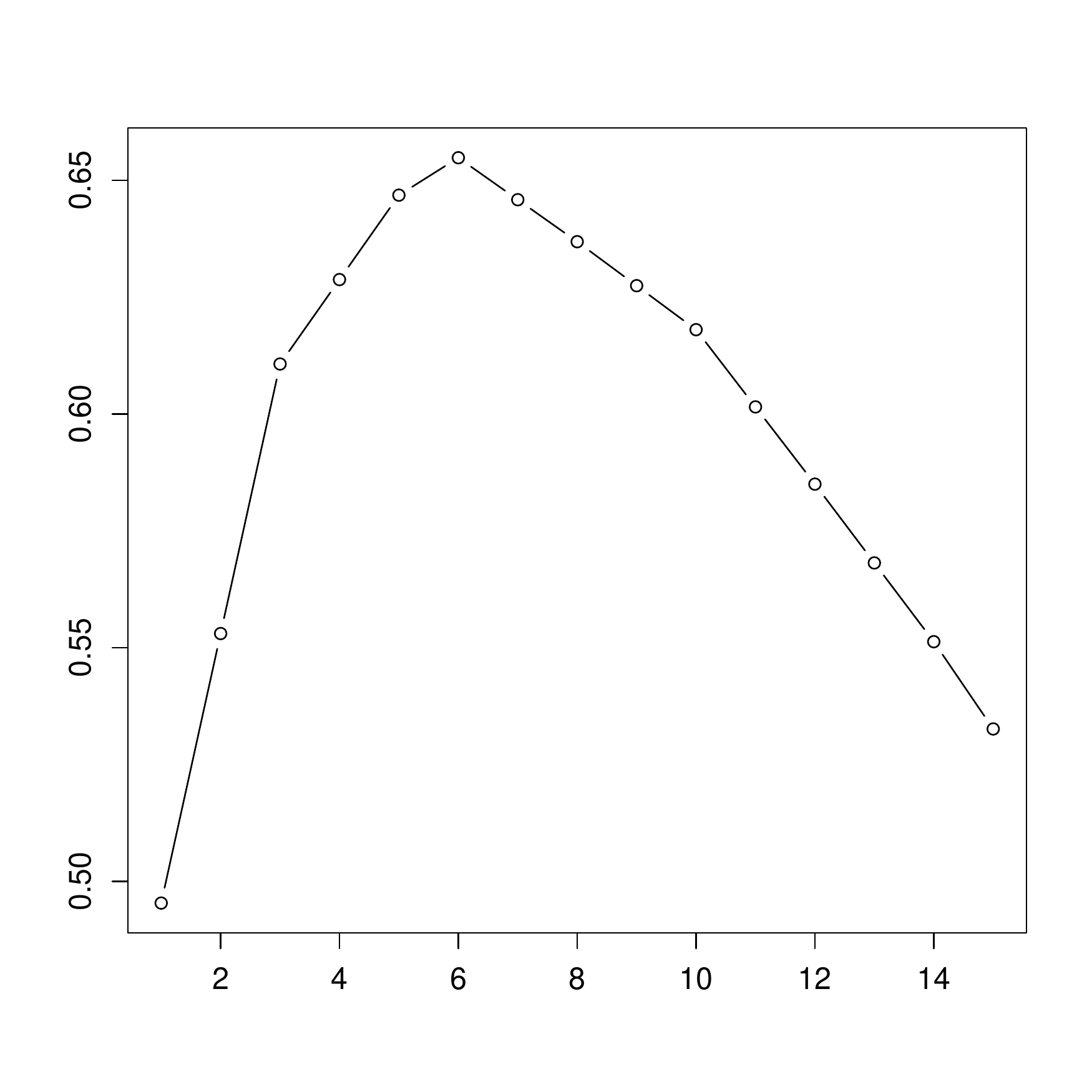}
 \caption{Adaptive filtering in the LP domain for the bipartite graph using the model selection scheme described in Section \ref{sec:filtering}.}
   \label{fig:bic}
\end{figure}
\subsection{Smooth Estimation by Adaptive Filtering} \label{sec:filtering}
The goal is to select the `significant' LP-transform coefficients $\LP[j,k;\cG]$ in a data-driven way. The data driven method of selecting the number of components is backed by substantial literature, see \cite{hart1992,ledwina94} and \cite{ledwina97}. Here, we apply the \cite{ledwina94} scheme using Schwarz selection criterion, which is known to be (1) consistent against \emph{every} alternative (including higher-frequency alternatives) to uniformity \citep{kallenberg95} and (2) enjoys neat asymptotic optimality property \citep{inglot2000}.

To identify the important elements using Ledwina's data-driven model selection criterion, we first rank order the squared LP-coefficients. Then we take the penalized cumulative sum of $k$ coefficients using information criterion $N^{-1}\log(N)k$ and choose the $k$ (number of LP components) for which this is maximum. Fig \ref{fig:bic} displays the penalized ordered cumulative sum of LP coefficients for bipartite graph, which selects the number of components $k=6$.

We estimate the LPINFOR statistic \eqref{eq:lpinfor} by taking the sum of squares of the selected LP-coefficients only. The number of filtered LP transforms (which denote the interesting graph correlations) determine the \emph{effective dimension} of the graph correlation density and tackle the spurious effects to produce smooth estimate $\widehat {\cd}(u,v;\cG)$.

\section{Applications}
\subsection{Diagnostic for Testing Null Model}
The model building process starts by checking whether the given process comes from the null model. In time series analysis, the flat shape of power spectral density and the autocorrelation function plot provide a quick and elegant diagnostic tool for the white noise process (which is the null model in the context of time series) that is widely practiced. The following question arises:

\vskip.15em
\emph{Can we propose an analogous portmanteau diagnostic test for null random graph models? What is a reasonable choice of null model for graph data?}
\vskip.15em
In the same spirit as time series analysis, here we propose two (graphical and nonparametric) diagnostic tests based on LP-Fourier transforms and graph correlation density function.

\begin{thm}[Asymptotic Normality]\label{thm:an} The sample LP graph transforms are i.i.d and
\beq \sqrt{N} \,\Big[\tLP(j,k;\cG) \,-\, \LP(j,k;\cG) \Big]\, \xrightarrow{d}\, \cN(0,1).\eeq
\end{thm}

Define the bivariate score function $L_{jk}:=S_j(u)S_k(v):[0,1]^2\rightarrow \cR$ and verify $\tLP(j,k;\cG)$ can be written as
$\Ex\big[L_{jk} \big( \tilde F(X;X), \tilde F(Y;Y) \big); \tilde p \big]$. The asymptotic normality result then follows immediately by using the same reasoning as given in \cite{ruymgaart1972}.

As a corollary to the previous theorem consider the following interesting fact: The graph LP transforms for random graph $\cG$ with probability of an edge equal to $p(x;X)p(y;Y)$ behaves asymptotically normal distribution with mean $0$ and variance $1/N$. This result has a threefold significance. First, it tells us the (theoretical) graph correlation density \eqref{eq:GCOR} is essentially ``flat'' on the unit square (compare with the power spectral density analogy of time series). Second, it suggests a reasonable definition for null models for graph data, which surprisingly matches the recommendation given in Section 3 of \cite{newman2006}. Third, there is a striking similarity between the asymptotic normality results for sample auto-correlations of white noise process (null model for time series) and our result. One difference is that our higher-order LP graph correlations have the flexibility to detect nonlinear graph connection patterns whereas auto-correlation is limited to linear patterns.

\begin{figure}[thb]
 \centering
 \vspace{-3em}
 \includegraphics[height=\textheight,width=.45\textwidth,keepaspectratio,trim=1cm 0cm 1cm 0cm]{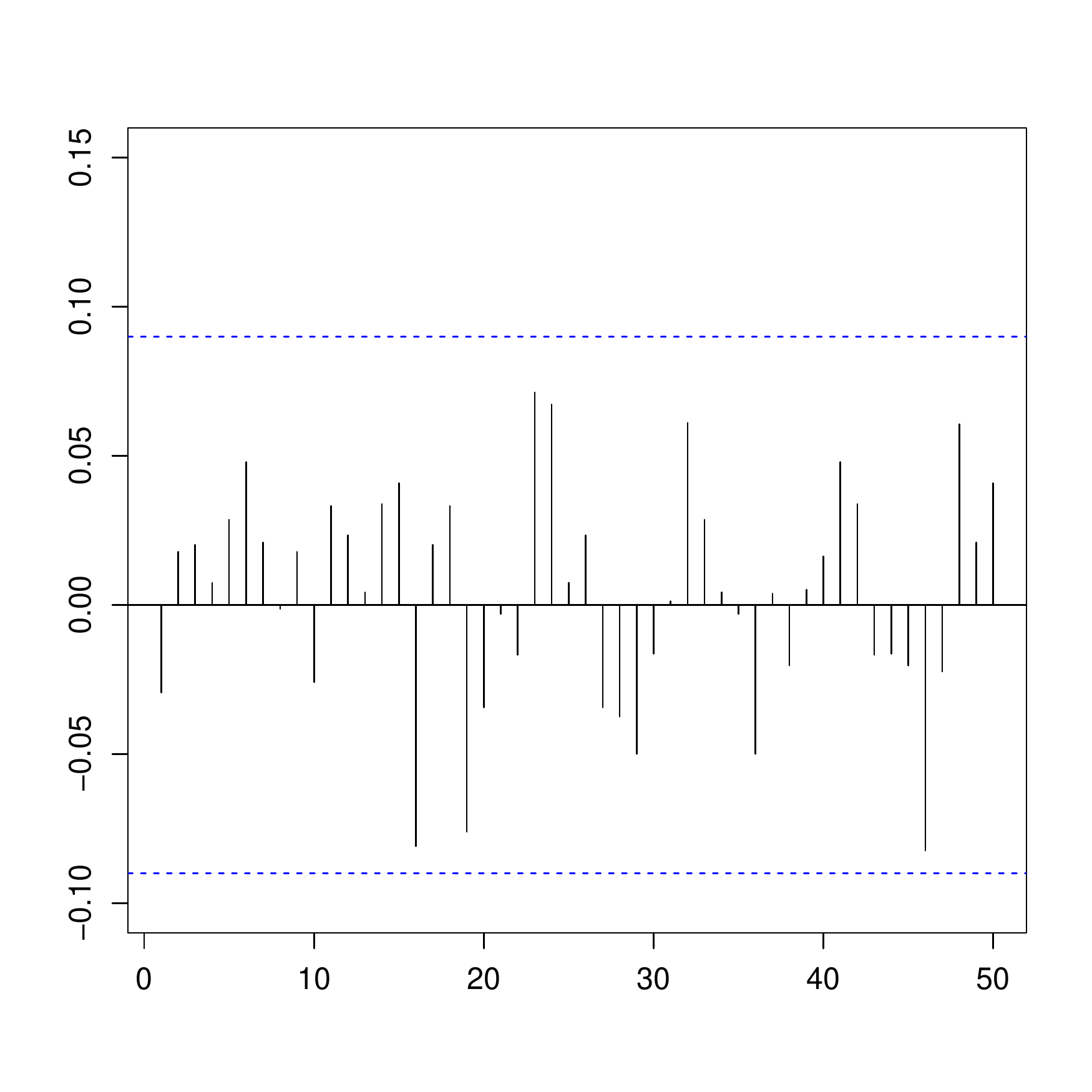}
 \includegraphics[height=\textheight,width=.54\textwidth,keepaspectratio,trim=.7cm 1cm 1.4cm 0cm]{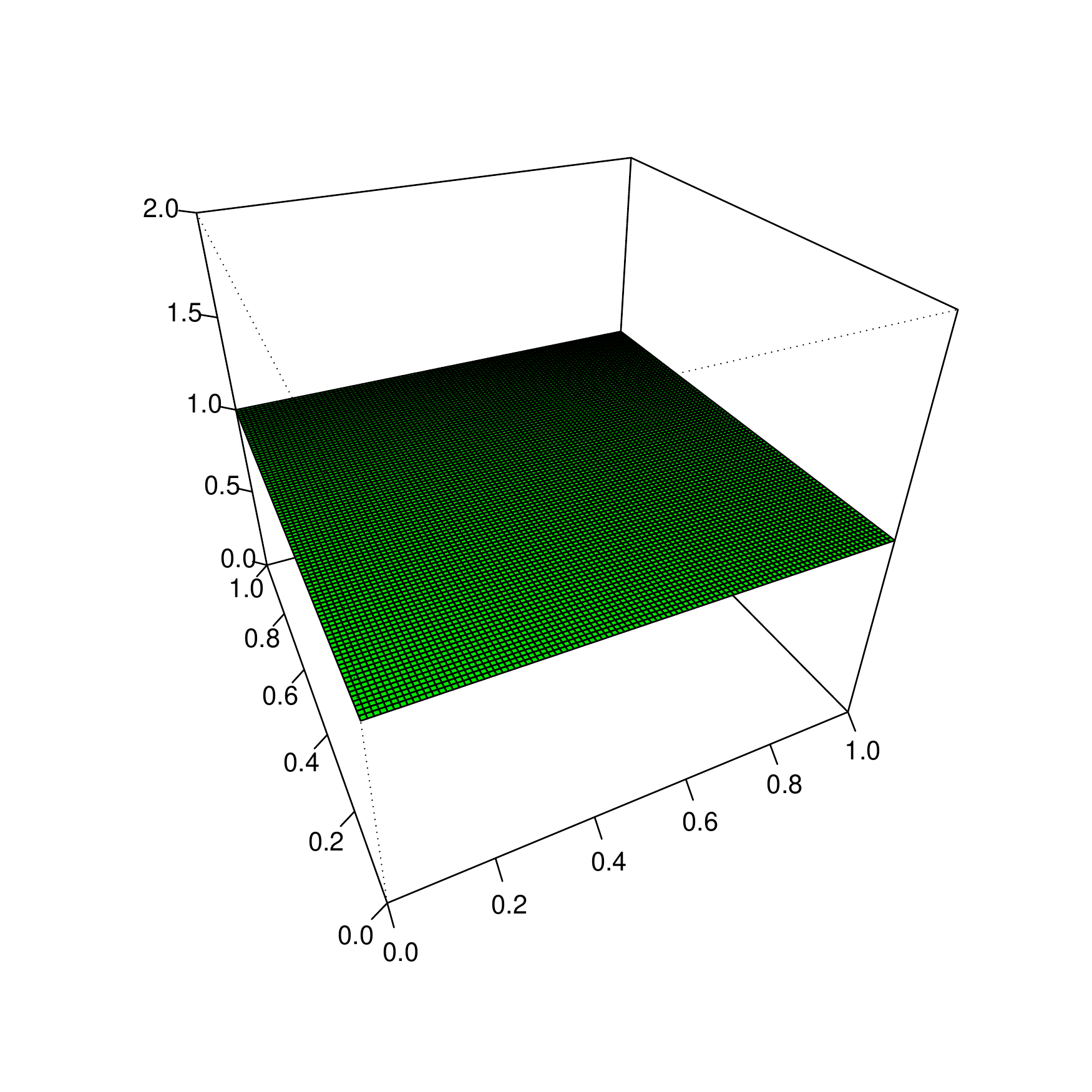}
 \caption{(color online) Graphical diagnostic for testing \ER graph (one specific null model); parallel to white noise testing in time series analysis. (a) LP-graph correlogram; (b) Smooth nonparametrically estimated graph correlation density.}
  \label{fig:testing}
   \vspace{-.5em}
\end{figure}
Figure \ref{fig:testing} plots the sample LP-transforms or LP-graph correlogram (along with the band $\pm 1.96/\sqrt{N} $, the $95\%$ confidence band) and the \emph{smooth} nonparametrically estimated graph correlation density for the \ER random graph model \citep{erdds1959}, a member of the null model class. The ``flat'' or uniform shape of the graph correlation density for \ER model signifies that the \emph{underlying strength of connectedness between all pairs of nodes are the same} (analogous to that in time series the flat spectrum of white noise indicates an equal amount of variance contained in all frequencies).
Various other diagnostic statistics for checking random graphs can be developed as a functional of $\cd(u,v;\cG)$; for example, consider the omnibus test statistic  $\LPINFOR(\cG)$, which can be shown to follow chi-square limit distribution under null using Theorem \ref{thm:an}. It is interesting to note that, our LPINFOR can be viewed as a generalization of \cite{hong1996} and Box and Pierce's (1970)\nocite{box1970} type test statistic for graph data.

\begin{figure}[ht]
 \centering
 \vspace{-3em}
 \includegraphics[height=\textheight,width=.45\textwidth,keepaspectratio,trim=1cm 0cm 1cm 0cm]{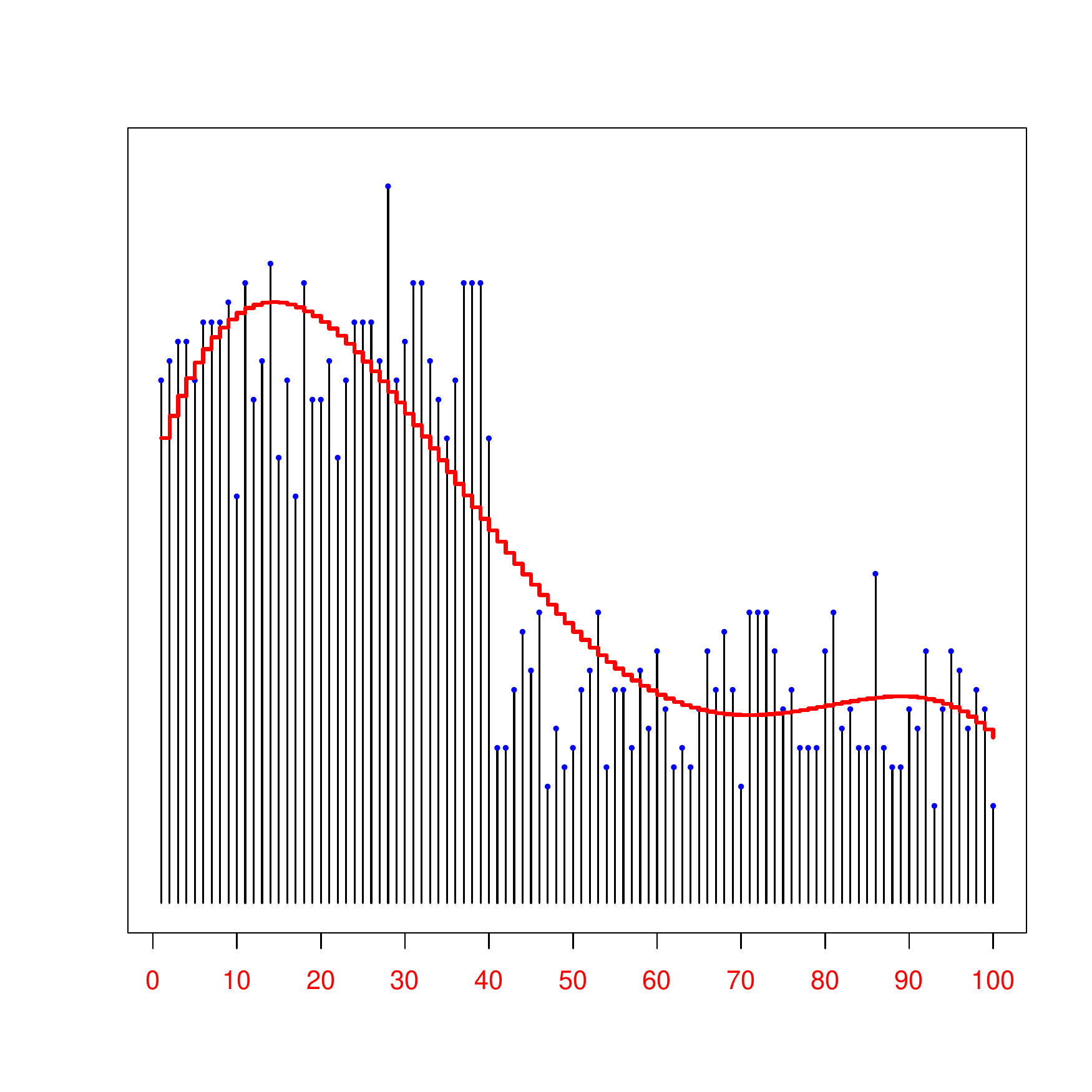}
 \includegraphics[height=\textheight,width=.54\textwidth,keepaspectratio,trim=.7cm 1.5cm 1.4cm 0cm]{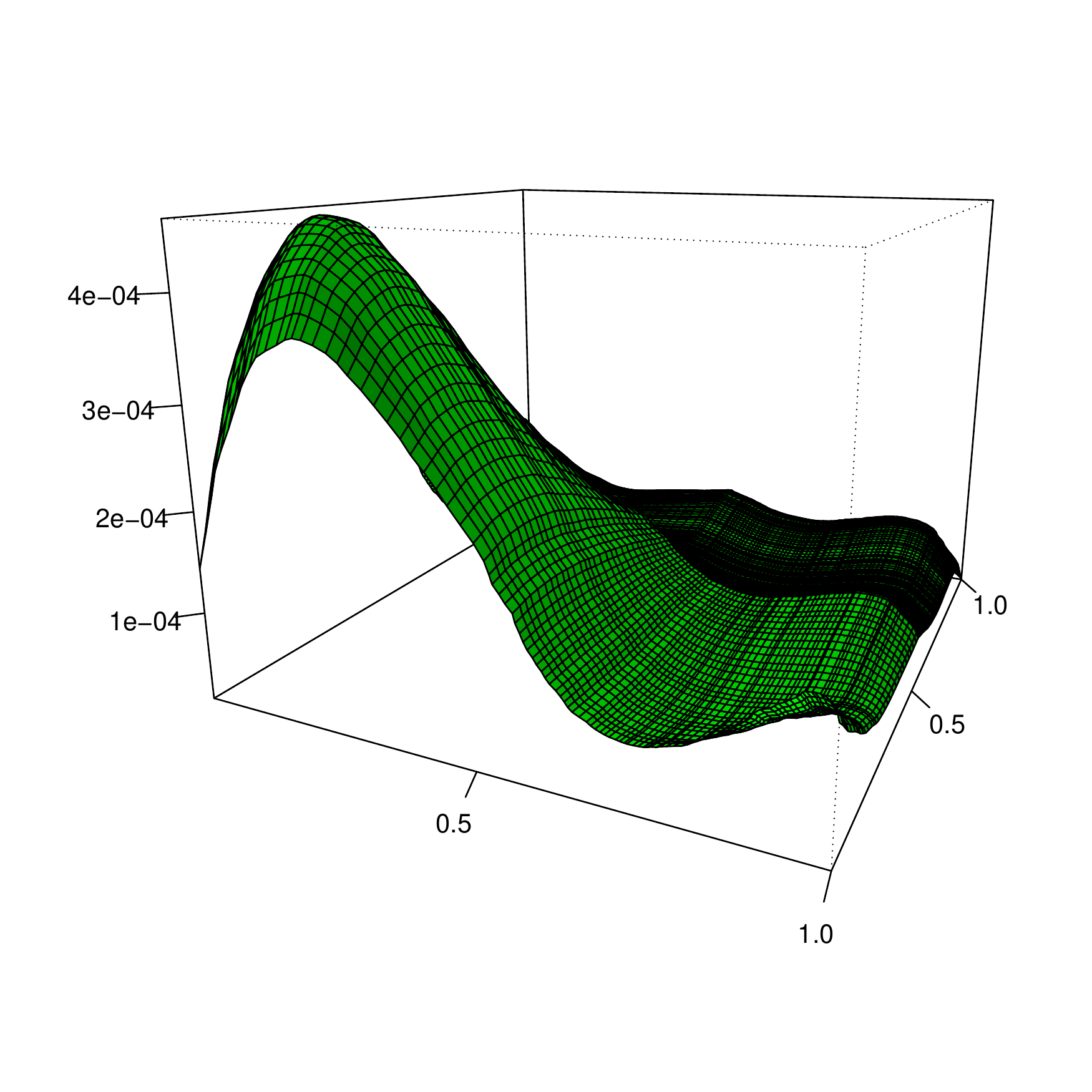}
 \caption{(color online) (a)left: smooth marginal estimate $\widehat{p}(x;X)$; (b)right: smooth nonparametric graphon estimate $\widehat{W}(u,v;\cG)$ of two-block stochastic model (with only one tightly connected group) with sizes $40$ and $60$ with the edge probability matrix $P=\big[ .6,.1\,;\,.1,.1 \big]$}
  \label{fig:graphon1}
\end{figure}
\subsection{Smooth Graphon Estimation} 
Nonparametric modeling of random graphs via Graphon is an emerging frontier of modern random graph theory. We define Graphon \citep{lovasz2012} as a bivariate function, $W(u,v;\cG):[0,1]^2 \rightarrow [0,1]$ denotes the probability of a directed edge from node $x$ to $y$, where $Q(u;X)=x$ and $Q(v;Y)=y$. The nonparametric estimation of this theoretical graph object was first initiated by \cite{olhede2013}, and \cite{airoldi14}, where they proposed histogram estimates by fitting a stochastic block model (to get the histogram bins by block-averaging). Unfortunately, this neither produces a smooth estimate nor a parsimonious modeling strategy. More specifically, we are interested in the following question:
\vskip.15em
\emph{How can we develop a nonparametric graphon approximation algorithm that will offer smooth estimate via sparse or low-rank modeling?}
\vskip.15em
Our LP estimation scheme makes the first attempt to do so by exploiting the novel connection between graphon and graph correlation density. The observation that every graphon can be written as a product of null graph model $p(x;X)p(y;Y)$ and the graph correlation density function $\cd(F(x;X),F(y;Y);\cG)$ for $x,y \in V(\cG)$ permits a new way to estimate graphon, which is applicable for directed/undirected graphs. Besides using nonparametrically estimated  $\widehat{\cd}(u,v;\cG)$, one may wish to use smooth estimates of the marginals (instead of raw empirical distributions) $\widehat{p}(x;X)$ and $\widehat{p}(y;Y)$ (expressed as linear combination of discrete LP polynomials as described in the next example).

Consider a two-block stochastic model (with only one tightly connected group) of sizes $40$ and $60$ with the edge probability matrix $P=\big[ .6,.1\,;\,.1,.1 \big]$. Use plug-in smooth nonparametric estimate of the graph correlation density to estimate the graphon. For a more enhanced answer, further smooth the marginals using the model: $p(x;X)=p_0(x;X) \big[ 1 + \sum_{j=1}^{n-1}\LP[j;p_0,p] T_j(x;p_0)\big]$, where $p_0(x;X)=n^{-1}$ is a discrete uniform distribution on the integers $x=\{1,2,\ldots,n\}$. Our specially design orthonormal polynomials $T_j(x;p_0)$ associated with discrete uniform distribution  $p_0(x;X)$ exactly match the $\DLeg_j$ discrete Legendre polynomials on $(0,1)$ (which belong to the `LP' family of polynomials). Fig \ref{fig:graphon1}(a) shows the fitted (smooth) nonparametric model $\widehat{p}(x;X)$ by selecting the components exactly the same way as discussed in Section 2.4. The resulting LP smooth graphon estimate is shown in Figure \ref{fig:graphon1}(b).
\begin{figure}[bht]
 \centering
 \vspace{-1em}
  \includegraphics[height=\textheight,width=.59\textwidth,keepaspectratio,trim=1.5cm 1.5cm 1.5cm 1.5cm]{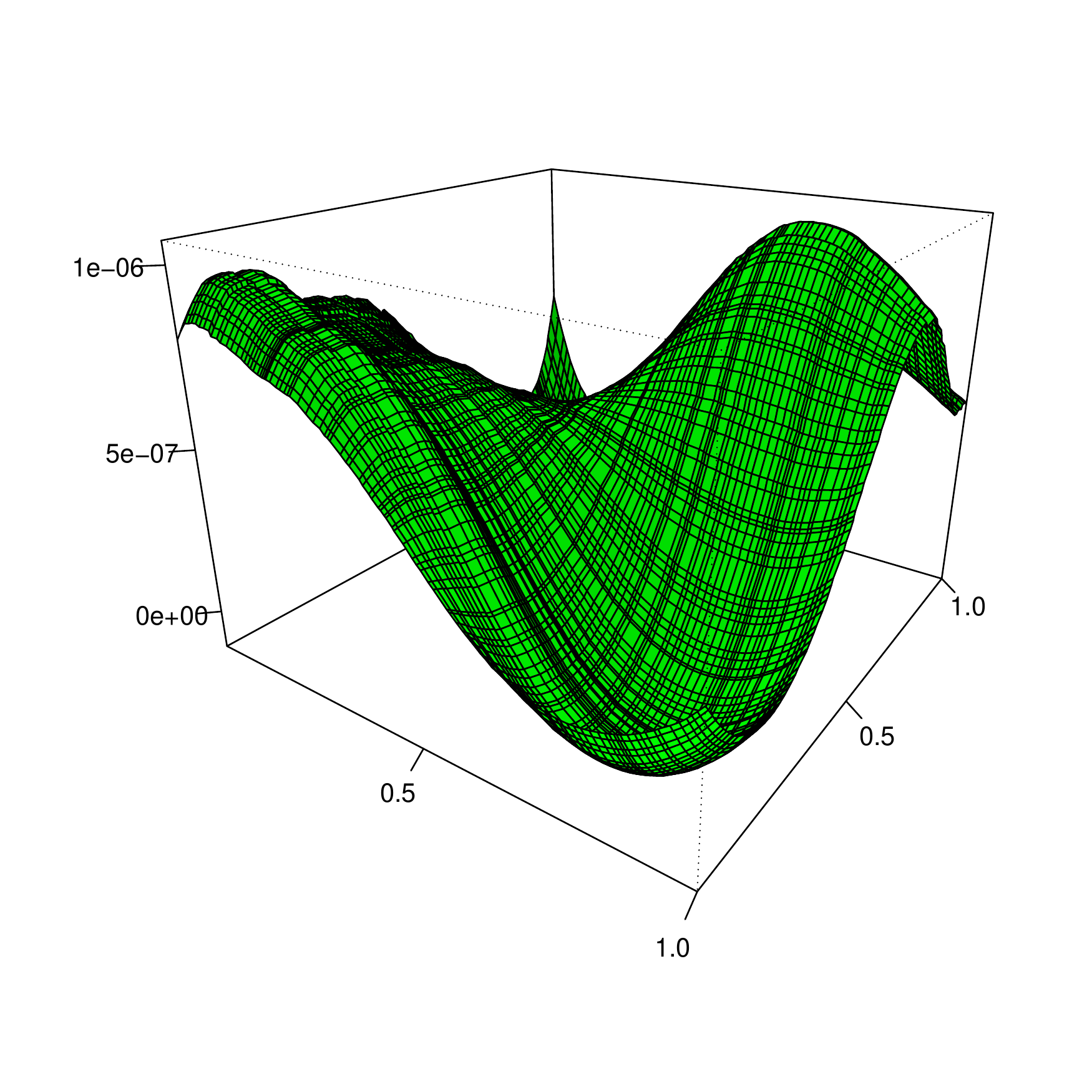}
     \vspace{-.25em}
 \caption{(color online) Smooth nonparametric graphon estimate of the political blog data, which is a directed sparse network.}
  \label{fig:graphon2}
   \vspace{-.25em}
\end{figure}

Figure \ref{fig:graphon2} shows the estimated asymmetric smooth graphon $\widehat{W}(u,v;\cG)$ for the political blog data. Bi-modality denotes the presence of two communities. One notable aspect of the shape of the estimated graphon is the contrasting sharpness of the two modes: one is ``flat'' and the other is more ``peaked''; this indicates the heterogeneity of strength of connections within each of the liberal and conservative political web-blogs.
%
%

\section{Conclusion}
A modern nonparametric modeling approach of graph data is reported here, which has its root in \cite{Deep14LP}. More specifically, this article introduces a new function, called Graph Correlation Density Field, that characterizes the `strength of connectedness'  between pairs of vertices of a random graph. Our motivation stems from the desire to develop the concepts of Fourier and inverse-Fourier like transformations for graph data to allow frequency domain interpretation that is similar to time series analysis. To accomplish this, we construct a new class of (data-driven) discrete orthonormal polynomials for graph data. We also discuss the issue of graph compression and sparse representation, which are especially useful for large graphs.  We show how these new concepts and tools can be used to solve many outstanding graph modeling problems. In particular, we focus on two problems in this article: (1) developing graphical diagnostic tests for null graph models in the same spirit of testing white noise process in time series analysis. Two diagnostic tools are introduced, which are analogous to the sample auto-correlation function plot and nonparametric spectral density based white noise testing approach for time series data; (2) the problem of smooth estimation of graphon is discussed by connecting it with graph correlation density. Throughout this paper, we pay special attention to develop ideas and algorithms that are analogous to the time series theory.

We hope this new perspective will enable us to develop a systematic and comprehensive theory of random graphs that is parallel to that of frequency-domain time series theory. The next article of this series will establish some surprising connections between our (modern nonparametric statistics based) LP graph theory and (matrix theory and linear algebra based) spectral graph theory.

\section*{Acknowledgment}
I would like to express my gratitude to Professor Emanuel Parzen for the encouragement to pursue this research direction. The author also thanks Mohsen Pourahmadi for constructive comments and several useful suggestions.
%

\vskip1em
\bib

\end{document}